\documentclass[11pt, reqno]{amsart}
\usepackage{graphicx,indentfirst}
\usepackage{amscd}

\usepackage[all]{xy}
\usepackage{amsmath}
\usepackage{amsfonts,amssymb}
\usepackage{amsthm}
\usepackage{latexsym}
\usepackage[american]{babel}
\usepackage{times}
\usepackage{tikz}
\usetikzlibrary{matrix,arrows,decorations.pathmorphing}
\usepackage[a4paper,left=2.5cm,right=2.5cm,bottom=3cm,top=3.5cm]{geometry}
\RequirePackage[T1]{fontenc}

\begin{document}
\title{Covariant Weil algebras}
\author{Zhaoting Wei}
\address{Department of Mathematics, 209 S 33 Street, University of Pennsylvania, Philadelphia, PA, 19104}
\email{zhaotwei@sas.upenn.edu}

\newcommand{\End}{\text{End}}
\newcommand{\ad}{\text{ad}}
\newcommand{\Ad}{\text{Ad}}
\newcommand{\Pol}{\text{Pol}}
\newcommand{\Cl}{\text{Cl}}
\newcommand{\B}{\text{B}}

\newtheorem{thm}{Theorem}[section]
\newtheorem{lemma}[thm]{Lemma}
\newtheorem{prop}[thm]{Proposition}
\newtheorem{coro}[thm]{Corollary}
\newtheorem{defi}{Definition}[section]
{\theoremstyle{plain} \newtheorem{rmk}{Remark}}

\maketitle

\tableofcontents

\begin{abstract}
 In this paper we introduce the classical and quantum covariant Weil algebras $W_{\tau}(\mathfrak{g})$ and $\mathcal{W}_{\tau}(\mathfrak{g})$. Covariant Weil algebras are simultaneous generalizations of Weil algebras (\cite{AM1}) and family algebras (\cite{Ki1}). We will define differentials, Lie derivatives and contractions on $W_{\tau}(\mathfrak{g})$ and $\mathcal{W}_{\tau}(\mathfrak{g})$ to make them curved-dg algebras. Moreover, the expression of curvatures will be given. It is hoped that covariant Weil algebras can be used in the construction of Mackey's analogue in \cite{Hig2}.
\end{abstract}

\section{Introduction}

The Weil algebra
$$
W(\mathfrak{g}^*):=S\mathfrak{g}^*\otimes \wedge \mathfrak{g}^*
$$ is an algebraic model of the differential forms on the universal principal bundle $\Omega^*(EG)$. Here $\mathfrak{g}$ is the Lie algebra of the Lie group $G$.

$W(\mathfrak{g}^*)$ is a $\mathfrak{g}$-differential graded algebra, that is , on $W(\mathfrak{g}^*)$ there exist the Lie derivative $L_a$, the contraction $\iota_a$ and the differential $d$, which satisfies certain relations, especially $d\circ d=0$.

In \cite{AM1}, A. Alekseev and E. Meinrenken introduce the noncommutative Weil algebra
$$
\mathcal{W}(\mathfrak{g}):=U(\mathfrak{g})\otimes \Cl(\mathfrak{g}),
$$
which is a deformation of $W(\mathfrak{g}^*)$. They shown that $\mathcal{W}(\mathfrak{g})$ is also a $\mathfrak{g}$-differential graded algebra.\\

In this paper, for a finite dimensional representation $\tau$ of $g$, we define the classical covariant Weil algebra and the quantum covariant Weil algebra
\begin{align*}
W_{\tau}(\mathfrak{g}^*):=&S\mathfrak{g}^*\otimes \wedge \mathfrak{g}^*\otimes \End V_{\tau},\\
\mathcal{W}_{\tau}(\mathfrak{g}):=&U(\mathfrak{g})\otimes\Cl(\mathfrak{g})\otimes \End V_{\tau}.
\end{align*}

$W_{\tau}(\mathfrak{g}^*)$ and $\mathcal{W}_{\tau}(g)$ are no longer $g$-differential graded algebra. Instead, they are curved $g$-differential graded algebra. That means they also have the Lie derivative, the contraction and the differential. However, $d\circ d$ is not $0$.

In fact we can find the curvature elements $C$ on $W_{\tau}(\mathfrak{g}^*)$ and $\mathcal{C}$ on $\mathcal{W}_{\tau}(\mathfrak{g})$ such that
$$d\circ d(-)=[C,-]$$ on $W_{\tau}(\mathfrak{g}^*)$; and $$d\circ d(-)=[\mathcal{C},-]$$ on $\mathcal{W}_{\tau}(\mathfrak{g})$.

The explicit formulas for $C$ and $\mathcal{C}$ will be given. Some of the properties will be studied.

\begin{rmk}
$W_{\tau}(\mathfrak{g}^*)$ can be considered as an algebraic model of homomorphism of vector bundles. The differential $d$ on $W_{\tau}(\mathfrak{g}^*)$ is corresponding to the covariant derivative on a vector bundle.
\end{rmk}

\subsection*{Acknowledgement}
The author would like to thank his advisor Jonathan Block for introducing me to curved-dg algebras, for his encouragement and helpful discussions. The author am also very grateful to Nigel Higson teaching me Mackey's analogue. Eckhard Meinrenken has explain to the author about the ideas of  noncommutative Weil algebra and the proof of Duflo's isomorphism theorem, which is very helpful. I would like to thank Alberto Garc\'{i}a-Raboso and Eric Korman for their helpful comments about this work.

\section{A review of the Weil algebra}\label{A review of the Weil algebra}
All the result in this section can be found in \cite{AM1}, section 2.

Let $\mathfrak{g}$ be a finite dimensional real algebra and $e_1,\ldots e_n$~be a basis. Let $f_{ij}^k$ be the structure constant, i.e
$$
[e_i, e_j]=f_{ij}^k e_k.
$$

Let $\mathfrak{g}^*$ be the dual linear space of $g$. The Weil algebra is defined to be
\begin{equation}
W(\mathfrak{g}^*):=S\mathfrak{g}^*\otimes \wedge \mathfrak{g}^*.
\end{equation}

$W(\mathfrak{g}^*)$ is a graded commutative super-algebra with the l-th grade
$$
W(\mathfrak{g}^*)^l=\bigoplus_{j+2k=l}S^k \mathfrak{g}^* \otimes \wedge^j \mathfrak{g}^*.
$$

Let  $e^i$~be a dual basis on  $\mathfrak{g}^*$. As in \cite{AM1}, we denote the corresponding element in  $S\mathfrak{g}^*$~by  $v^i$,~and the corresponding element in  $\wedge \mathfrak{g}^*$~by  $y^i$.

There are natural Lie derivative on  $S\mathfrak{g}^*$~and  $\wedge \mathfrak{g}^*$.~On generators
\begin{equation}
\begin{split}
L_a y^c=&-f^c_{ab}y^b,\\
L_a v^c=&-f^c_{ab}v^b.
\end{split}
\end{equation}

On the whole  $W(\mathfrak{g}^*)$~we define the Lie derivative as:
\begin{equation}
L_a:=L_a\otimes 1+1\otimes L_a.
\end{equation}

$L_a$~is a derivation of degree  $0$.

There is also a natural contraction on  $\wedge \mathfrak{g}^*$.~On $W(\mathfrak{g}^*)$~we define the contraction by
\begin{equation}
\iota_a:=1\otimes \iota_a.
\end{equation}

$\iota_a$~is a derivation of degree  $-1$.

Finally we have the differential. On generators:
\begin{equation}
\begin{split}
d^W y^a:=& v^a-\frac{1}{2}f^a_{jk}y^jy^k,\\
d^W v^a:=& -f^a_{jk}y^jv^k.
\end{split}
\end{equation}

$d^W$ is a derivation of degree  $1$.

It is easy to check that  $d^W\circ d^W=0,~[L_a, d^W]=0,~[L_a, \iota_b]=f^c_{ab}\iota_c$ and the Cartan's formula  $[\iota_a, d^W]=L_a$ holds. Notice that
 $[\iota_a, d^W]=\iota_a d^W+d^W\iota_a$ is the super-commutator.

The weil algebra  $W(\mathfrak{g}^*)$ is an algebraic model of the universal principal bundle in topology and geometry.

\section{The classical covariant Weil algebra}
\subsection{The definition of  the classic covariant Weil algebra $W_{\tau}(\mathfrak{g}^*)$}
\begin{defi}\label{classic covariant weil algebra}let  $(\tau, V_{\tau})$~be a finite dimensional representation of  $g$.~Let us define the classical covariant Weil algebra  $W_{\tau}(\mathfrak{g}^*)$~associated with  $\tau$~to be:
\begin{equation}
W_{\tau}(\mathfrak{g}^*):=S\mathfrak{g}^*\otimes \wedge \mathfrak{g}^*\otimes \End V_{\tau}.
\end{equation}

In the following content we use  $\tau_a$~to denote $\tau(e_a)\in \End V_{\tau}$~and  $A$~to denote a general element in $\End V_{\tau}$.

$W_{\tau}(\mathfrak{g}^*)$~is a graded algebra with elements in  $\End V_{\tau}$~of degree  $0$.
\end{defi}

\begin{rmk}
Unlike $W(\mathfrak{g})$, $W_{\tau}(\mathfrak{g}^*)$ is not graded commutative.
\end{rmk}

\subsection{The three operators on  $W_{\tau}(\mathfrak{g}^*)$}\label{The three operators on W}
\begin{defi}\label{The Lie derivative on  W tau g*}
The Lie derivative on  $W_{\tau}(\mathfrak{g}^*)$~ is defined to be
\begin{equation}
L_a:=L_a\otimes 1\otimes 1+1\otimes L_a\otimes 1+1\otimes 1\otimes \Ad \tau_a.
\end{equation}
$L_a$~is a derivation of degree  $0$.
\end{defi}

\begin{defi}\label{The contraction on  W tau g*}
The contraction only acts on the second component:
\begin{equation}
\iota_a:=1\otimes \iota_a\otimes 1.
\end{equation}
$\iota_a$~is a derivation of degree  $-1$.
\end{defi}

Now we come to the step to define the covariant differential on  $W_{\tau}(\mathfrak{g}^*)$. We define  $d^{W,\tau} y^a$~and  $d^{W,\tau} v^a$~exactly the same as in  $W(\mathfrak{g}^*)$.~As for  $A\in \End(V_{\tau})$,~we define
$$
d^{W,\tau} A:= 1\otimes y^a\otimes [\tau_a, A].
$$

In conclusion, we have
\begin{defi}\label{The covariant differential on  W tau g*}
The covariant differential $d^{W,\tau}$ on  $W_{\tau}(\mathfrak{g}^*)$ has the following expression
\begin{equation}
\begin{split}
d^{W,\tau} y^a:=& v^a-\frac{1}{2}f^a_{jk}y^jy^k,\\
d^{W,\tau} v^a:=& -f^a_{jk}y^jv^k,\\
d^{W,\tau} A:=& 1\otimes y^a\otimes [\tau_a, A].
\end{split}
\end{equation}
$d^{W,\tau}$~is a derivation of degree  $1$.
\end{defi}

\begin{rmk}
This definition is inspired by the differential in the Cartan-Eilenberg algebra.
\end{rmk}

\begin{rmk}
It is clear that if we restrict to  $W(\mathfrak{g}^*)$,~the $L_a,\iota_a$~and $d^{W,\tau}$ coincide with the corresponding original operators .
\end{rmk}

Next we can compute the commutators of the three operators. We can see that  $[L_a, d^{W,\tau}]=0$ still holds on  $W_{\tau}(\mathfrak{g}^*)$. It is sufficient to check it on the generators. For $y^a$~and $v^a$~the equality is the same as for  $W(\mathfrak{g}^*)$. For $A\in \End(V_{\tau})$:
\begin{equation}
\begin{split}
[L_a, d^{W,\tau}]A=&L_a d^{W,\tau} A-d^{W,\tau}L_a A\\
=& L_a(y^b[\tau_b,A])-d^{W,\tau}([\tau_a, A])\\
=& L_a(y^b)[\tau_b,A]+y^b(L_a[\tau_b,A])-y^c[\tau_c,[\tau_a, A]]\\
=& -f^b_{ac}y^c[\tau_b,A]+y^b[\tau_a,[\tau_b,A]]-y^c[\tau_c,[\tau_a, A]]\\
=&-f^b_{ac}y^c[\tau_b,A]+y^b[\tau_a,[\tau_b,A]]-y^c[\tau_a,[\tau_c,A]]+y^c[[\tau_a,\tau_c], A]~(\text{Jacobi identity})
\\
=&-f^b_{ac}y^c[\tau_b,A]+y^b[\tau_a,[\tau_b,A]]-y^c[\tau_a,[\tau_c,A]]+f^b_{ac}y^c[\tau_b,A]\\
=&0
\end{split}
\end{equation}

By definition we get $[L_a, \iota_b]=f^c_{ab}\iota_c$.

Finally let us check the Cartan's formula  $[\iota_a, d^{W,\tau}]=L_a$~on  $W_{\tau}(\mathfrak{g}^*)$~. For $y^a$~and $v^a$~the equality is the same as for  $W(\mathfrak{g}^*)$. For $A$~
\begin{equation}
\begin{split}
[\iota_a, d^{W,\tau}]A=&\iota_a d^{W,\tau} A+d^{W,\tau}\iota_a A\\
=&\iota_a d^{W,\tau} A\\
=&\iota_a(y^k[\tau_k, A]) \text{ (\,We omit the tensor symbols.) }\\
=&[\tau_a, A]\\
=& L_a A.
\end{split}
\end{equation}

\subsection{$d^{W,\tau}\circ d^{W,\tau}$~and the curvature of  $W_{\tau}(\mathfrak{g}^*)$}

Now the main problem to us is to compute $d^{W,\tau}\circ d^{W,\tau}$. On  $y^a$~and $v^a$~, $d^{W,\tau}\circ d^{W,\tau}=0$~holds as it holds in  $W(\mathfrak{g}^*)$.
However, on  $A\in \End(V_{\tau})$~we have:
\begin{equation}\label{dd for A classic}
\begin{split}
d^{W,\tau}\circ d^{W,\tau}A=&d^{W,\tau}(y^a[\tau_a, A])\\
=&d^{W,\tau} (y^a)[\tau_a, A]-y^a d^{W,\tau}([\tau_a, A])\\
=& (v^a-\frac{1}{2}f^a_{jk}y^jy^k)[\tau_a, A]-y^ay^b[\tau_b,[\tau_a, A]]\\
=& v^a[\tau_a, A]-\frac{1}{2}f^a_{jk}y^jy^k[\tau_a, A]-y^ay^b[\tau_b,[\tau_a, A]]\\
=& v^a[\tau_a, A]-\frac{1}{2}f^a_{jk}y^jy^k[\tau_a, A]+\frac{1}{2}y^ay^b[[\tau_a,\tau_b],A]~(\text{By Jacobi identity})\\
=&v^a[\tau_a, A]-\frac{1}{2}f^a_{jk}y^jy^k[\tau_a, A]+\frac{1}{2}f^c_{ab}y^ay^b[\tau_c, A]\\
=&v^a[\tau_a, A]
\end{split}
\end{equation}

\begin{defi}[The curvature of  $W_{\tau}(\mathfrak{g}^*)$]\label{curvature for cov weil classic}
Let $C:=v^a\tau_a\in W_{\tau}(\mathfrak{g}^*)$.~It is obvious that  $C$~is independent of the choice of the basis of  $g$. $C$~is called the \textbf{curvature} of  $W_{\tau}(\mathfrak{g}^*)$~
 \end{defi}

 We have the following:

\begin{prop}\label{dd for W(g, tau) classic}
$d^{W,\tau}\circ d^{W,\tau}(-)=[C,-]$~on  $W_{\tau}(\mathfrak{g}^*)$.
\end{prop}
Proof:It is sufficient to check it on the generators. On $y^a$~and $v^a$~both sides are $0$. On $A\in \End(V_{\tau})$~it has been done in ~(\ref{dd for A classic}).
$\square$
\\

The curvature  $C$~has the following important property:
\begin{prop}\label{bianchi for W(g, tau) classic}
$C$~is closed in  $W_{\tau}(\mathfrak{g}^*)$, i.e.  $d^{W,\tau} C=0$.
\end{prop}
Proof: Direct computation shows:
\begin{equation}
\begin{split}
d^{W,\tau} C=&d^{W,\tau} (v^a\tau_a)\\
=& d^{W,\tau}(v^a)\tau_a+v^a d^{W,\tau}(\tau_a)\\
=& -f^a_{jk}y^jv^k\tau_a+ v^ay^b[\tau_b, \tau_a]\\
=& -f^a_{jk}y^jv^k\tau_a+f^c_{ba}v^ay^b\tau_c\\
=& 0~(\text{Notice that}~y^j~\text{commutes with}~v^k) ~\square
\end{split}
\end{equation}

\begin{rmk}Proposition \ref{bianchi for W(g, tau) classic} is the algebraic version of the Bianchi identity in differential geometry.
\end{rmk}

Now we see that  $W_{\tau}(\mathfrak{g}^*)$~ is a curved differential graded algebra with $\mathfrak{g}$ action on it (\,for short, curved $\mathfrak{g}$-dga).

\subsection{The flat elements in $W_{\tau}(\mathfrak{g}^*)$}
In the study of curved differential graded algebras, it is important to find out which element of  $W_{\tau}(\mathfrak{g}^*)$ is  $0$~under  $d^{W,\tau}\circ d^{W,\tau}$. By proposition~\ref{dd for W(g, tau) classic}, $d^{W,\tau}\circ d^{W,\tau}(-)=[C,-]$. Therefore, the problem is equivalent to find the commutant of  $C$~in  $W_{\tau}(\mathfrak{g}^*)$.

\begin{defi}\label{flat in W(g, tau)}
Let
\begin{equation}
W_{\tau}(\mathfrak{g}^*)_F:=\{x\in W_{\tau}(\mathfrak{g}^*)|~d^{W,\tau}d^{W,\tau} x=[C,x]=0\}
\end{equation}
be the elements which commute with $C$ in $W_{\tau}(\mathfrak{g}^*)$. We also call them the flat elements in $W_{\tau}(\mathfrak{g}^*)$.
\end{defi}

We have the following result on $W_{\tau}(\mathfrak{g}^*)_F$:
\begin{prop}\label{sub curved g dga classic}
$W_{\tau}(\mathfrak{g}^*)_F$ is a curved sub $g$-differential graded algebra of~$W_{\tau}(\mathfrak{g}^*)$,~i.e. $W_{\tau}(\mathfrak{g}^*)_F$ is closed under addition, multiplication and the actions of $L_a$, $\iota_a$ and $d^{W,\tau}$.
\end{prop}
Proof: $W_{\tau}(\mathfrak{g}^*)_F$ is obviously closed under addition and $d^{W,\tau}$. Since $[C,-]$ is a derivation, $W_{\tau}(\mathfrak{g}^*)_F$, as the commutant of $C$, is also closed under multiplication.

Since $[L_a, d^{W,\tau}]=0$ on $W_{\tau}(\mathfrak{g}^*)$, i.e. $L_a d^{W,\tau}=d^{W,\tau}L_a$, it is easy to see that $W_{\tau}(\mathfrak{g}^*)_F$ is closed under $L_a$.

Lastly, since $[\iota_a, d^{W,\tau}]=L_a$, $\forall x\in W_{\tau}(\mathfrak{g}^*)_F$, i.e. $d^{W,\tau}d^{W,\tau}x=0$, we have
\begin{align*}
d^{W,\tau}d^{W,\tau}(\iota_a x)=&d^{W,\tau}(L_a x-\iota_a d^{W,\tau}x)\\
=& d^{W,\tau}L_a x-d^{W,\tau}\iota_adx\\
=& d^{W,\tau}L_a x-(L_a-\iota_ad^{W,\tau})d^{W,\tau}x\\
=& d^{W,\tau}L_ax-L_ad^{W,\tau}x+\iota_a d^{W,\tau}d^{W,\tau}x\\
=& [d^{W,\tau}, L_a]x+\iota_a(d^{W,\tau}d^{W,\tau}x)\\
=&0
\end{align*}

Hence $\iota_a x$ is also in $W_{\tau}(\mathfrak{g}^*)_F$, i.e. $W_{\tau}(\mathfrak{g}^*)_F$ is closed under $\iota_a$. $\square$
\\

We have proved in Propostion \ref{bianchi for W(g, tau) classic} that $C$ itself is in $W_{\tau}(\mathfrak{g}^*)_F$. It is also obvious that all the scalars (more precisely, scalar matrices) are in $W_{\tau}(\mathfrak{g}^*)_F$. Since $y^j$ commutes with $u^k$ and $\tau_a$, we also have $\wedge \mathfrak{g}^*\subset W_{\tau}(\mathfrak{g}^*)_F$.

In fact we can show more.

\begin{defi}\label{Horizontal  and basic calssic Weil algebra and commutant}
Let
\begin{equation}
W_{\tau}(\mathfrak{g}^*)_{hor}:=S\mathfrak{g}^*\otimes \End V_{\tau}
\end{equation}
 be the horizontal elements in $W_{\tau}(\mathfrak{g}^*)$,
 \begin{equation}
W_{\tau}(\mathfrak{g}^*)_{basic}:=\{x\in W_{\tau}(\mathfrak{g}^*)_{hor}|L_a x=0\}
\end{equation}
  be the basic elements in $W_{\tau}(\mathfrak{g}^*)$, and
 \begin{equation}
W_{\tau}(\mathfrak{g}^*)_{hor,F}:=\{x\in W_{\tau}(\mathfrak{g}^*)_{hor}|~[C,x]=0\}
\end{equation}
be the flat elements in $W_{\tau}(\mathfrak{g}^*)$.
\end{defi}

By definition, $W_{\tau}(\mathfrak{g}^*)_{hor,F} \subset W_{\tau}(\mathfrak{g}^*)_F$ and $W_{\tau}(\mathfrak{g}^*)=\wedge \mathfrak{g}^*\otimes W_{\tau}(\mathfrak{g}^*)_{hor}$. Furthermore we can prove the following:

\begin{prop}\label{decompose of W(g, tau)_F}
$W_{\tau}(\mathfrak{g}^*)_F=\wedge \mathfrak{g}^*\otimes W_{\tau}(\mathfrak{g}^*)_{hor,F}$
\end{prop}
Proof: Since $W_{\tau}(\mathfrak{g}^*)=\wedge \mathfrak{g}^*\otimes W_{\tau}(\mathfrak{g}^*)_{hor}$, for any $x\in W_{\tau}(\mathfrak{g}^*)_F$, we can write
$$
x=\sum_i f_is_i, \text{ here } f_i\in \wedge \mathfrak{g}^*, s_i\in W_{\tau}(\mathfrak{g}^*)_{hor}.
$$

We combine the terms so that the $f_i$'s  are linear independent.

Since $f_i\in \wedge \mathfrak{g}^*$, we get $[C, f_i]=0$. Hence we have
$$
0=[C,x]
=[C,\sum_i f_is_i]
=\sum_i f_i[C,s_i]
$$

Since the $f_i$'s are linear independent as we arranged before, we get $[C,s_i]=0$ for all $i$.

Hence $W_{\tau}(\mathfrak{g}^*)_F=\wedge \mathfrak{g}^*\otimes W_{\tau}(\mathfrak{g}^*)_{hor,F}$. $\square$\\

By proposition~\ref{decompose of W(g, tau)_F}, to study $W_{\tau}(\mathfrak{g}^*)_F$ it is sufficient to study $W_{\tau}(\mathfrak{g}^*)_{hor,F}$.

Obviously, $S\mathfrak{g}^*\in W_{\tau}(\mathfrak{g}^*)_{hor,F}$. In fact we can say more:

\begin{prop}\label{classic basic is flat}
We have
$W_{\tau}(\mathfrak{g}^*)_{basic}\subset W_{\tau}(\mathfrak{g}^*)_{hor,F}$
\end{prop}
Proof: Let $f^i\otimes \theta_i\in W_{\tau}(\mathfrak{g}^*)_{basic}$, hence
$$
f^i\otimes[\tau_a,\theta_i]=-L_a f^i\otimes \theta_i \text{ for each } a.
$$
For $C=u^a\tau_a$, we have
\begin{align*}
[C, f^i\theta_i]=&[u^a\tau_a,f^i\theta_i]\\
=&u^af^i[\tau_a,\theta_i]\\
=&u^a(f^i[\tau_a,\theta_i])\\
=&-u^aL_af^i\theta_i
\end{align*}

Now the proposition reduces to the

\begin{lemma}\label{classic casimir}
Let $u^a, L_a$ are defined as usual, then $\forall f\in S\mathfrak{g}^*$, we have
\begin{equation}
u^aL_af=0
\end{equation}
\end{lemma}
Proof of the lemma \ref{classic casimir}: It is obvious that $u^aL_a$ is a derivation on $S\mathfrak{g}^*$, so it is sufficient to
prove the identity for the generators of $S\mathfrak{g}^*$. We see that
$$
u^aL_a u^b=u^a(-f^b_{ac}u^c)=-f^b_{ac}u^au^c.
$$
Since $f^b_{ac}$ is anti-symmetric with respect to $a,c$, the sum is $0$. $\square$

By lemma \ref{classic casimir}, we have $[C, f^i\theta_i]=0$, i.e. $[C, f^i\theta_i]\in W_{\tau}(\mathfrak{g}^*)_{hor,F}$. $\square$

\begin{rmk}
For quadratic Lie algebra $g$, $W_{\tau}(\mathfrak{g}^*)_{basic}$ is the same as the \textbf{classical family algebra} $\mathcal{C}_{\tau}(\mathfrak{g})$ introduced by A.A. Kirillov (see \cite{Ki1}) and our proposition \ref{classic basic is flat} is equivalent to the lemma 1 in section 1 of \cite{Ki1}.
\end{rmk}

Now we have
\begin{equation}
S\mathfrak{g}^*\cdot W_{\tau}(\mathfrak{g}^*)_{basic}\in W_{\tau}(\mathfrak{g}^*)_{hor,F}.
\end{equation}
We don't know whether $S\mathfrak{g}^*\cdot W_{\tau}(\mathfrak{g}^*)_{basic}= W_{\tau}(\mathfrak{g}^*)_{hor,F}$. For the case $g=sl(2, \mathbb{R})$ and $\tau$ be the standard or adjoint representations of $sl(2, \mathbb{R})$, it is true by the structure theorem of family algebras (see \cite{Roz} chapter 6 ). But it is still unknown in the general cases.

\section{A review of the noncommutative Weil algebra}\label{A review of the noncommutative (quantum) Weil algebra}
All the result in this section can be found in \cite{AM1} section 3.

First we introduce the quadratic Lie algebra
\begin{defi}\label{quadratic lie algebra}
A Lie algebra $\mathfrak{g}$ is called a quadratic Lie algebra if there exists a non-degenerate symmetric bilinear form $\B$ which is $\mathfrak{g}$-invariant, i.e.
\begin{equation}
\B([x,y],z)+\B(y,[x,z])=0, \forall x,y,z\in \mathfrak{g}.
\end{equation}
\end{defi}

It is well-known that a semisimple Lie algebra is quadratic.\\

Let $e_a$ be a basis for $g$, $[e_a,e_b]=f^c_{ab}e_c$ as usual. $\B(e_a, e_b)=\B_{ab}$, Then $B_{ab}=B_{ba}$
and
\begin{equation}\label{g-invariance of bilinear form}
f^c_{ab}B_{cd}+f^c_{ad}B_{bc}=0
\end{equation}

We can use $\B$ to lower the index of the structure constant: $f_{abc}:=\B_{ad}f^d_{bc}$.\\

For a quadratic Lie algebra we can indentify $g$ and $\mathfrak{g}^*$ via $\B$, so
$$
S\mathfrak{g}^*\otimes \wedge \mathfrak{g}^*\cong S\mathfrak{g}\otimes \wedge \mathfrak{g}
$$

If we treat $e^a$ as a basis of $g$ instead of $\mathfrak{g}^*$, then the $L_a$, $\iota_a$, $d$ can be defined on $Sg\otimes \wedge g$ using the same formula as on $S\mathfrak{g}^*\otimes \wedge \mathfrak{g}^*$. Let
$$
W(\mathfrak{g}):=S\mathfrak{g}\otimes \wedge \mathfrak{g}.
$$
We have $W(\mathfrak{g})\cong W(\mathfrak{g}^*)$ as a $\mathfrak{g}$-dga.\\

For a quadratic Lie algebra $\mathfrak{g}$ we can construct a noncommutative analogue of $W(\mathfrak{g})$. First we introduce the Clifford algebra:

\begin{defi}\label{Clifford algebra}
The Clifford algebra associated to $\mathfrak{g}$ is defined to be
\begin{equation}
\Cl (\mathfrak{g}):=T(\mathfrak{g})/<\B(\mu,\nu)-\mu\otimes \nu-\nu\otimes \mu>.
\end{equation}
\end{defi}

Now we come to the noncommutative (quantum) Weil algebra is

\begin{defi}\label{noncommutative Weil algebra}
\begin{equation}
\mathcal{W}(\mathfrak{g}):=U(\mathfrak{g})\otimes \Cl(\mathfrak{g}).
\end{equation}
\end{defi}

\begin{rmk}
$\mathcal{W}(\mathfrak{g})$ is called noncommutative Weil algebra in \cite{AM1} and \cite{AM2}. In this paper I prefer to call it quantum Weil algebra to distinguish $\mathcal{W}(\mathfrak{g})$ from $W_{\tau}(\mathfrak{g})$, since $W_{\tau}(\mathfrak{g})$ is also noncommutative.
\end{rmk}

$\mathcal{W}(\mathfrak{g})$ is a filtered algebra where
$$
\mathcal{W}^{(l)}(\mathfrak{g}):= \sum_{2j+k=l}U^{(j)}(\mathfrak{g})\otimes \Cl^{(k)}(\mathfrak{g})
$$
and $\mathcal{W}(\mathfrak{g})$ has a well defined $\mathbb{Z}_2$-grade. So we can talk about even and odd elements in $\mathcal{W}(\mathfrak{g})$.\\

For $\mathcal{W}(\mathfrak{g})$ we also have the three operators $L_a, \iota_a$ and $d$ as for $W(\mathfrak{g})$.

\begin{rmk} A significant difference between $\mathcal{W}(\mathfrak{g})$ and $W(\mathfrak{g})$ is that the three operators on $\mathcal{W}(\mathfrak{g})$ are all inner derivations.
\end{rmk}

Follow \cite{AM1}, we choose an orthonormal basis under $\B$ to simplify the computation, i.e.
\begin{equation}
\\B(e_a,e_b)=\delta_{ab}.
\end{equation}
Hence $f_{abc}=f^a_{bc}$.

Now let $u_a$ denote the $e_a$ in $U(\mathfrak{g})$ and $x_a$ denote the $e_a$ in $\Cl(\mathfrak{g})$. Define (Notice that we take sum for small Latin letter indices.)
\begin{equation}
\begin{split}
g_a:=&-\frac{1}{2}f_{ars}x_rx_s\in \Cl^{(2)}(\mathfrak{g}).\\
\gamma:=&\frac{1}{3}x_ag_a=-\frac{1}{6}f_{abc}x_ax_bx_c\in \Cl^{(3)}(\mathfrak{g}).\\
\mathfrak{D}:=& x_au_a+\gamma\in \mathcal{W}^{(3)}(\mathfrak{g})
\end{split}
\end{equation}

Then we come to the

\begin{defi}\label{three operator for noncommutative weil algebra}
On $\mathcal{W}(\mathfrak{g})$ we have the Lie derivative, the contraction and the differential map:
\begin{equation}
\begin{split}
L_a:=&\ad(u_a+g_a)\\
\iota_a:=&\ad(x_a)\\
d^{\mathcal{W}}:=&\ad(\mathfrak{D}).
\end{split}
\end{equation}

It is easy to see that
\begin{equation}
\begin{split}
L_a x_b=&f_{cab}x_c,~L_a u_b=f_{cab}u_c,\\
\iota_a x_b=& \delta_{ab},~\iota_a u_b=0,\\
d^{\mathcal{W}}x_a=&u_a-\frac{1}{2}f_{abc}x_bx_c,\\
d^{\mathcal{W}}u_a=&-f_{abc}x_bu_c.
\end{split}
\end{equation}

So these definitions make sense and are analogous to those on $W(\mathfrak{g})$.
\end{defi}

On $\mathcal{W}(\mathfrak{g})$ we also have the relations $d^{\mathcal{W}}\circ d^{\mathcal{W}}=0,~[L_a, d^{\mathcal{W}}]=0,~[L_a, \iota_b]=f^c_{ab}\iota_c$~and the Cartan's formula  $[\iota_a, d^{\mathcal{W}}]=L_a$.

In fact, by the (super) Jacobi identity, the relations above reduce to the corresponding commutation relations between $x_a, g_a, u_a, \gamma$ and $\mathfrak{D}$. For those commutators we have :
\begin{lemma}\label{commutators in quantum W(g)}
\begin{equation}
\begin{split}
&[x_a, g_b]=-[g_b, x_a]=-f_{bac}x_c\\
&[u_a+g_a,\mathfrak{D}]=0\\
&[x_a,\mathfrak{D}]=u_a+g_a.
\end{split}
\end{equation}
\end{lemma}
Proof: See \cite{AM1} section 3. $\square$\\

Here we give the detail for $d^{\mathcal{W}}\circ d^{\mathcal{W}}=0$ to illustrate the story, although the computation can be found in \cite{AM1} section 3. The proof is as follows:
\begin{equation}
\begin{split}
d^{\mathcal{W}}\circ d^{\mathcal{W}}(-)=&[\mathfrak{D},[\mathfrak{D},-]]\\
=&[\frac{1}{2}[[\mathfrak{D},\mathfrak{D}],-]~\text{(by super) Jacobi identity}\\
=&[\mathfrak{D}^2,-].
\end{split}
\end{equation}

So it is sufficient to compute $\mathfrak{D}^2=\frac{1}{2}[[\mathfrak{D},\mathfrak{D}]$ in $\mathcal{W}(\mathfrak{g})$.
Remember that we take sum for small Lation letter indices.
\begin{equation}
\begin{split}
\mathfrak{D}^2=\frac{1}{2}[[\mathfrak{D},\mathfrak{D}]=&\frac{1}{2}[u_ax_a+\gamma,u_ax_a+\gamma]\\
=&\frac{1}{2}[u_ax_a,u_ax_a]+\frac{1}{2}[u_ax_a,\gamma]+\frac{1}{2}[\gamma,u_ax_a]+\frac{1}{2}[\gamma,\gamma]\\
=&\frac{1}{2}[u_ax_a,u_bx_b]+[u_ax_a,\gamma]+\frac{1}{2}[\gamma,\gamma]\\
=&\frac{1}{2}[u_a,u_b]x_ax_b+\frac{1}{2}u_au_b[x_a,x_b]+u_a[x_a,\gamma]+\gamma^2\\
=& \frac{1}{2}f_{cab}u_cx_ax_b+\frac{1}{2}u_au_b\delta_{ab}+u_a[x_a, \gamma]+\gamma^2\\
=&-u_cg_c+\frac{1}{2}u_au_a+u_a[x_a, \gamma]+\gamma^2.
\end{split}
\end{equation}

By computation as in  \cite{AM1} section 3, we have the following:
\begin{lemma}\label{gamma xa and gamma square}
\begin{equation}
\begin{split}
&[x_a, \gamma]=g_a\\
&\gamma^2= -\frac{1}{48}f_{abc}f_{abc}~~\text{is a scalar}
\end{split}
\end{equation}
\end{lemma}
Proof of the lemma: In \cite{AM1} section 3. $\square~$\\

So
\begin{equation}
\mathfrak{D}^2=\frac{1}{2}u_au_a+\gamma^2
\end{equation}

Now we notice that $\frac{1}{2}u_au_a$ is the Casimir~~element, which belongs to the center of $U(\mathfrak{g})$, hence
$\mathcal{W}(\mathfrak{g})$. As a scalar, $\gamma^2$ also belongs to the center of $\mathcal{W}(\mathfrak{g})$. So
\begin{equation}
d^{\mathcal{W}}\circ d^{\mathcal{W}}(-)=[\frac{1}{2}u_au_a+\gamma^2,-]=0.
\end{equation}

\section{The quantum covariant Weil algebra}

\subsection{The definition of  $\mathcal{W}_{\tau}(\mathfrak{g})$}

\begin{defi}\label{quantum covariant weil algebra}
Let  $(\tau, V_{\tau})$~be a finite dimensional representation of  $g$. Similar to the commutative case, we
define the quantum covariant Weil algebra  $\mathcal{W}_{\tau}(\mathfrak{g})$~associated with  $\tau$~to be:
\begin{equation}
\mathcal{W}_{\tau}(\mathfrak{g}):=U(\mathfrak{g})\otimes\Cl(\mathfrak{g})\otimes \End V_{\tau}.
\end{equation}

Again we use  $\tau_a$~to denote $\tau(e_a)\in \End V_{\tau}$~and  $A$~to denote a general element in $\End V_{\tau}$.

$\mathcal{W}_{\tau}(\mathfrak{g})$ is a filtered and $\mathbb{Z}_2$-graded algebra with elements in  $\End V_{\tau}$ of degree  $0$.
\end{defi}

\subsection{The three operators on  $\mathcal{W}_{\tau}(\mathfrak{g})$}
Inspired by the constructions in section \ref{The three operators on W} and section \ref{A review of the noncommutative (quantum) Weil algebra}, we can construct the Lie derivative, contraction and covariant differential on $\mathcal{W}_{\tau}(\mathfrak{g})$.
\begin{defi}\label{three operators on quantum covariant weil algebra}
\begin{equation}
\begin{split}
L_a:=&\ad(u_a+g_a+\tau_a)\\
\iota_a:=&\ad(x_a)\\
d^{\mathcal{W},\tau}:=&\ad(\mathfrak{D}+x_a\tau_a)
\end{split}
\end{equation}
\end{defi}

It is easy to see that, when restricted on $\mathcal{W}(\mathfrak{g})\subset \mathcal{W}_{\tau}(\mathfrak{g})$, $L_a$ and $\iota_a$ are the same as the original $L_a$ and $\iota_a$ defined in section \ref{A review of the noncommutative (quantum) Weil algebra}.

On $A\in \End V_{\tau}$, $L_a A=[\tau_a, A]$ and $\iota_a A=0$ as we expect.

\begin{rmk}
$d^{\mathcal{W},\tau}$ is \textbf{not} equal to $d^{\mathcal{W}}$ even if we restrict to $\mathcal{W}(\mathfrak{g})$. This is a significant difference between classic and quantum covariant Weil algebras.
\end{rmk}

In fact, for $x_b$ we have
\begin{equation}
\begin{split}
d^{\mathcal{W},\tau}x_b=&[\mathfrak{D}+x_a\tau_a, x_b]\\
=&[\mathfrak{D}, x_b]+[x_a\tau_a,x_b]\\
=&d^{\mathcal{W}}x_b+[x_a,x_b]\tau_a\\
=&d^{\mathcal{W}}x_b+(\iota_a x_b)\tau_a
\end{split}
\end{equation}

Notice that on $u_b$, we have$d^{\mathcal{W},\tau}u_b=d^{\mathcal{W}}u_b$. So we conclude that on $\mathcal{W}(\mathfrak{g})\subset \mathcal{W}_{\tau}(\mathfrak{g})$,
\begin{equation}
d^{\mathcal{W},\tau}=d^{\mathcal{W}}+\iota_a\tau_a
\end{equation}

On $A\in \End V_{\tau}$, as we expect,
\begin{equation}
d^{\mathcal{W},\tau}A=x_a[\tau_a,A].
\end{equation}

\subsection{The commutation relation of $L_a$, $\iota_a$ and $d^{\mathcal{W},\tau}$}
Similar to the relations on $\mathcal{W}(\mathfrak{g})$, the commutators of $L_a$, $\iota_a$ and $d^{\mathcal{W},\tau}$ reduce to the commutators of $u_a+g_a+\tau_a$, $x_a$ and $\mathfrak{D}+x_a\tau_a$, by the super-Jacobi identity. It requires some careful calculation.

First of all
\begin{prop}
$$
[L_a,\iota_b]=f_{abc}\iota_c.
$$
\end{prop}
Proof:
\begin{equation}
\begin{split}
[L_a,\iota_b](-)=&[[u_a+g_a+\tau_a, x_b],-]\\
=&[[g_a,x_b],-]\\
=&[f_{abc}x_c, -] \text{ (\,Lemma } \ref{commutators in quantum W(g)})\\
=&f_{abc}\iota_c(-). ~\square
\end{split}
\end{equation}

Next we have
\begin{prop}
$$
[L_a,d^{\mathcal{W},\tau}]=0.
$$
\end{prop}
Proof: \begin{equation}
\begin{split}
[L_a,d^{\mathcal{W},\tau}](-)=&[[u_a+g_a+\tau_a, \mathfrak{D}+x_b\tau_b],-]\\
=&[[u_a+g_a+\tau_a,\mathfrak{D}]+[u_a+g_a+\tau_a,x_b\tau_b],-]
\end{split}
\end{equation}
As in Lemma \ref{commutators in quantum W(g)}, we have $[u_a+g_a,\mathfrak{D}]=0$ and obviously $[\tau_a, \mathfrak{D}]=0$, so
\begin{equation}
\begin{split}
\text{The above formula}=&[[u_a+g_a+\tau_a,x_b\tau_b],-]\\
=&[[g_a+\tau_a,x_b\tau_b],-]\\
=&[[g_a,x_b]\tau_b]+x_b[\tau_a,\tau_b],-]\\
=&[f_{abc}x_c\tau_b+f_{abc}x_b\tau_c,-]\\
=&0.
\end{split}
\end{equation}

So
$$
[L_a,d^{\mathcal{W},\tau}]=0. ~\square
$$

Finally,
\begin{prop} We have the Cartan's formula:
$$
[\iota_a,d^{\mathcal{W},\tau}]=L_a.
$$
\end{prop}
Proof:
\begin{equation}
\begin{split}
[\iota_a,d^{\mathcal{W},\tau}](-)=&[[x_a, \mathfrak{D}+x_b\tau_b],-]\\
=& [[x_a,\mathfrak{D}]+[x_a,x_b\tau_b],-]\\
=& [u_a+g_a+[x_a,x_b]\tau_b,-]~(\,\text{Lemma}~\ref{commutators in quantum W(g)})\\
=&[u_a+g_a+\delta_{ab}\tau_b,-]\\
=&[u_a+g_a+\tau_a,-]=L_a(-).
\end{split}
\end{equation}

So
$$
[\iota_a,d^{\mathcal{W},\tau}]=L_a. ~\square
$$

\subsection{$d^{\mathcal{W},\tau}\circ d^{\mathcal{W},\tau}$~and the curvature of  $\mathcal{W}_{\tau}(\mathfrak{g})$}
Similar to $\mathcal{W}(\mathfrak{g})$
\begin{equation}
d^{\mathcal{W},\tau}d^{\mathcal{W},\tau}(-)=[\frac{1}{2}[\mathfrak{D}+x_a\tau_a,\mathfrak{D}+x_b\tau_b],-]
\end{equation}
So we need to compute $[\mathfrak{D}+x_a\tau_a,\mathfrak{D}+x_b\tau_b]$.
$$
[\mathfrak{D}+x_a\tau_a,\mathfrak{D}+x_b\tau_b]=[\mathfrak{D},\mathfrak{D}]+[\mathfrak{D},x_b\tau_b]
+[x_a\tau_a,\mathfrak{D}]+[x_a\tau_a,x_b\tau_b]
$$

In section \ref{A review of the noncommutative (quantum) Weil algebra} we have shown that
$$
[\mathfrak{D},\mathfrak{D}]=u_au_a+2\gamma^2.
$$
We also have
\begin{align*}
[\mathfrak{D},x_b\tau_b]+[x_a\tau_a,\mathfrak{D}]=&2[x_a\tau_a,\mathfrak{D}]\\
=&2[x_a,\mathfrak{D}]\tau_a\\
=&2(u_a+g_a)\tau_a.
\end{align*}
The last term is
\begin{align*}
[x_a\tau_a,x_b\tau_b]=&[x_a,x_b]\tau_a\tau_b+x_ax_b[\tau_a,\tau_b]\\
=&\delta_{ab}\tau_a\tau_b+x_ax_bf_{cab}\tau_c\\
=&\tau_a\tau_a-2g_c\tau_c.
\end{align*}

Sum them up we get
\begin{prop}\label{dd for W(g, tau) quantum}
\begin{equation}
\begin{split}
[\mathfrak{D}+x_a\tau_a,\mathfrak{D}+x_b\tau_b]=&u_au_a+2\gamma^2+2(u_a+g_a)\tau_a+\tau_a\tau_a-2g_c\tau_c\\
=&u_au_a+2u_a\tau_a+\tau_a\tau_a+2\gamma^2.
\end{split}
\end{equation}

Therefore
\begin{equation}
d^{\mathcal{W},\tau}d^{\mathcal{W},\tau}(-)=[\frac{1}{2}(u_au_a+2u_a\tau_a+\tau_a\tau_a+2\gamma^2),-].
\end{equation}
\end{prop}
Proof: $\square$

\begin{defi}[The quantum curvature of  $\mathcal{W}_{\tau}(\mathfrak{g})$]\label{curvature for cov weil quantum}
Let
$$\mathcal{C}:=\frac{1}{2}(u_au_a+2u_a\tau_a+\tau_a\tau_a+2\gamma^2)\in \mathcal{W}_{\tau}(\mathfrak{g}).$$
It is obvious that  $\mathcal{C}$~is independent of the choice of the basis of  $g$. $\mathcal{C}$~is called the \textbf{quantum curvature} of  $\mathcal{W}_{\tau}(\mathfrak{g})$. We have proved that
$$
d^{\mathcal{W},\tau}d^{\mathcal{W},\tau}(-)=[\mathcal{C},-].
$$
 \end{defi}

Similar to the closednees of the curvature of $W_{\tau}(\mathfrak{g})$, we have the following result:
\begin{prop}\label{bianchi for W(g, tau) quantum}
$\mathcal{C}$~is closed in  $\mathcal{W}_{\tau}(\mathfrak{g})$, i.e.  $d^{\mathcal{W},\tau} \mathcal{C}=0$.
\end{prop}
Proof:
\begin{equation}
\begin{split}
d^{\mathcal{W},\tau} \mathcal{C}=&[\mathfrak{D}+x_a\tau_a,\mathcal{C}]\\
=&[\mathfrak{D},\frac{1}{2}[\mathfrak{D},\mathfrak{D}]]\\
=&\frac{1}{2}[\mathfrak{D},[\mathfrak{D},\mathfrak{D}]]
\end{split}
\end{equation}

By the super-Jacobi identity, we obtain $3[\mathfrak{D},[\mathfrak{D},\mathfrak{D}]]=0$ and hence $[\mathfrak{D},[\mathfrak{D},\mathfrak{D}]]=0$. As a result:
$$
d^{\mathcal{W},\tau} \mathcal{C}=0.~\square
$$

Same as  $W_{\tau}(\mathfrak{g}^*)$, $\mathcal{W}_{\tau}(\mathfrak{g})$ is a curved $\mathfrak{g}$-dga.

\subsection{The flat elements in $\mathcal{W}_{\tau}(\mathfrak{g})$}
Similar to the $W_{\tau}(\mathfrak{g})$ case, we have
\begin{defi}\label{flat in quantum W(g, tau)}
Let
\begin{equation}
\mathcal{W}_{\tau}(\mathfrak{g})_F:=\{x\in \mathcal{W}_{\tau}(\mathfrak{g})|~d^{\mathcal{W},\tau}d^{\mathcal{W},\tau} x=[\mathcal{C},x]=0\}
\end{equation}
be the elements which commute with $\mathcal{C}$ in $\mathcal{W}_{\tau}(\mathfrak{g}^*)$. We also call them the flat elements in $\mathcal{W}_{\tau}(\mathfrak{g}^*)$.
\end{defi}

We also have
\begin{prop}\label{sub curved g dga quantum}
$\mathcal{W}_{\tau}(\mathfrak{g})_F$~is a curved sub $\mathfrak{g}$-differential graded algebra of $\mathcal{W}_{\tau}(\mathfrak{g})$, i.e. it is closed under addition, multiplication and the actions of $L_a$, $\iota_a$ and $d^{\mathcal{W},\tau}$.
\end{prop}
Proof: The same as Proposition \ref{sub curved g dga classic}. $\square$

\begin{defi}\label{Horizontal  and basic quantum Weil algebra and commutant}
Let
\begin{equation}
 \mathcal{W}_{\tau}(\mathfrak{g})_{hor}:=U(\mathfrak{g})\otimes \End V_{\tau}
\end{equation}
 be the horizontal elements in $\mathcal{W}_{\tau}(\mathfrak{g}^*)$,
\begin{equation}
\mathcal{W}_{\tau}(\mathfrak{g})_{basic}:=\{x\in \mathcal{W}_{\tau}(\mathfrak{g})_{hor}|L_a x=0\}
\end{equation}
  be the basic elements in $\mathcal{W}_{\tau}(\mathfrak{g})$, and
\begin{equation}
\mathcal{W}_{\tau}(\mathfrak{g})_{hor,F}:=\{x\in \mathcal{W}_{\tau}(\mathfrak{g})_{hor}|~[\mathcal{C},x]=0\}
\end{equation}
be the flat elements in $\mathcal{W}_{\tau}(\mathfrak{g})$.
\end{defi}

\begin{rmk}$\mathcal{W}_{\tau}(\mathfrak{g})_{basic}$ is the same as the \textbf{quantum family algebra} $\mathcal{Q}_{\tau}(\mathfrak{g})$ introduced by A.A. Kirillov (see \cite{Ki1})
\end{rmk}

Similar to Proposition \ref{decompose of W(g, tau)_F}, we have the following

\begin{prop}\label{decompose of quantum W(g, tau)_F}
$\mathcal{W}_{\tau}(\mathfrak{g})_F=\Cl(\mathfrak{g})\otimes \mathcal{W}_{\tau}(\mathfrak{g})_{hor,F}$
\end{prop}
Proof: Since the elements in $\Cl(\mathfrak{g})$ commute with $\mathcal{C}$, the proof is the same as that of Proposition \ref{decompose of W(g, tau)_F}. $\square$

As a result, to study $\mathcal{W}_{\tau}(\mathfrak{g})_F$ it is sufficient to study $\mathcal{W}_{\tau}(\mathfrak{g})_{hor,F}$.

In fact we know very few about $\mathcal{W}_{\tau}(\mathfrak{g})_{hor,F}$. It is obvious that $U(\mathfrak{g})\nsubseteq \mathcal{W}_{\tau}(\mathfrak{g})_{hor,F}$. I don't know whether or not $\mathcal{W}_{\tau}(\mathfrak{g})_{basic}\subset \mathcal{W}_{\tau}(\mathfrak{g})_{hor,F}$.

\section{Further topics}

Much about covariant Weil algebras are still unknown. Here I list some of them.

\subsection{A suitable cohomology theory}
A differential graded algebra has its cohomology, which is naturally an algebra. Now our $W_{\tau}(\mathfrak{g})$
and $\mathcal{W}_{\tau}(\mathfrak{g})$ are curved-dga's, so they don't have cohomologies in the naive sense.

 Now the question is: Can we find a theory which substitute the cohomology theory and compatible with the $\mathfrak{g}$-action, i.e. Lie derivative and contraction, in our case?

\begin{rmk}
The theory of curved-dga is of increasing importance, for example see \cite{Block1} and \cite{Orlov1}. $W_{\tau}(\mathfrak{g})$
and $\mathcal{W}_{\tau}(\mathfrak{g})$ may provide good example to testify the theories.
\end{rmk}

\subsection{The quantization map between $W_{\tau}(\mathfrak{g})$
and $\mathcal{W}_{\tau}(\mathfrak{g})$}
In \cite{Hig2} Section 9, N. Higson proposed the problem of constructing a quantization map $Q$ between $\mathcal{C}_{\tau}(\mathfrak{g})$ and $\mathcal{Q}_{\tau}(\mathfrak{g})$ such that the following diagram commutes.
\begin{equation}
\begin{tikzpicture}[node distance=3cm, auto]
  \node (A1) {$\mathcal{C}_{\tau}(\mathfrak{g})$};
  \node [right of=A1] (A2) {$\mathcal{Q}_{\tau}(\mathfrak{g})$};
  \node [below of=A1, node distance=1.5cm] (B1) {$S(\mathfrak{h})$};
  \node [below of=A2, node distance=1.5cm] (B2) {$U(\mathfrak{h})$};

\draw[->] (A1) to node {$\text{GHC}_{\tau,c}$} (B1);
\draw[->,dashed] (A1) to node {$Q$}(A2);
\draw[->] (A2) to node {$\text{GHC}_{\tau}$}(B2);
\draw[->] (B1) to node {$\cong$}(B2);
\end{tikzpicture}
\end{equation}
Here $Q$ is a vector space isomorphism but need not to be an algebraic homomorphism.

It is a natural generalization of the famous Duflo's isomorphism theorem.

In \cite{AM1} and \cite{AM2} A. Alekseev, and E. Meinrenken give a new proof of Duflo's isomorphism theorem using the quantization map of the Weil algebras.

Since the covariant Weil algebras are generalization of the Weil algebras, it is expected that we can construct a quantization map between $W_{\tau}(\mathfrak{g})$ and $\mathcal{W}_{\tau}(\mathfrak{g})$ which gives the quantization map expected by Higson.

\end{document}